\documentclass[reprint,showkeys,amsmath,amssymb,aps]{revtex4-1}

\usepackage{comment}
\usepackage[usenames,dvipsnames]{xcolor}
\usepackage{graphicx}
\usepackage{cancel}
\usepackage{amsmath}
\usepackage{amssymb}
\usepackage{bm}

\usepackage{pgfplots}
\usetikzlibrary{plotmarks}

\begin{document}

\title{Stochastic switching in slow-fast systems: a large fluctuation approach}

\author{Christoffer~R.~Heckman}
\email{christoffer.heckman.ctr@nrl.navy.mil}

\author{Ira~B.~Schwartz}
\email{ira.schwartz@nrl.navy.mil}

\affiliation{U.S.\ Naval Research Laboratory, Code 6792\\
Plasma Physics Division, Nonlinear Dynamical Systems Section\\
Washington, DC 20375, USA}

\definecolor{dkgreen}{rgb}{0,0.6,0}
\definecolor{gray}{rgb}{0.5,0.5,0.5}
\definecolor{mauve}{rgb}{0.58,0,0.82}

\begin{abstract}
In this paper we develop a perturbation method to predict the rate of 
occurrence of rare events for singularly perturbed stochastic systems 
using a probability density function approach. In contrast to a stochastic
normal form approach, we model rare event 
occurrences due to large fluctuations probabilistically and 
employ a WKB ansatz to approximate their rate of occurrence. This results
in the generation of a two-point boundary value problem that models the 
interaction of the state variables and the most 
likely noise force required to induce a rare event. The resulting equations of 
motion of describing the phenomenon are shown to be 
singularly perturbed. Vastly different time scales among the 
variables are leveraged to reduce the dimension 
and predict the dynamics on the slow manifold in a deterministic 
setting. The resulting 
constrained equations of motion may be used to directly compute an
exponent that determines the probability of rare events.

To verify the theory, a stochastic damped Duffing oscillator with three 
equilibrium points (two sinks separated by a saddle) is analyzed. The 
predicted switching time between states is computed using 
the optimal path that resides in an expanded phase space. We show 
that the exponential scaling of the switching rate as a function of 
system parameters agrees well with numerical simulations. Moreover, 
the dynamics of the original system and the reduced system via center 
manifolds are shown to agree in an exponentially scaling sense.
\end{abstract}

\keywords{singular perturbation, stochastic differential equation, 
optimal path, noise, rare event}

\maketitle

\section{Introduction}
Many stochastic systems of physical interest possess dynamics
which occur over multiple time scales. These systems present unique
difficulties since the multiple time scales interact with the 
stochasticity to affect the dynamics, leading to phenomena such as 
stochastic switching resulting from large fluctuations. For 
deterministic slow-fast systems singular
perturbation theory may be applied to guide analysis, while noisy
systems are better understood through tools from statistical
mechanics.

The study of slow-fast systems has recently become popular as a 
result of the insight it affords into fields such as chemical 
reactions and electro-mechanical systems \cite{desroches}. Due to the 
presence of 
distinct timescales on which phenomena occur in these singularly 
perturbed systems, it becomes mathematically tractable to apply 
perturbation methods to accurately predict the behavior of high-order 
systems in terms of low order ones. This model reduction greatly
simplifies bifurcation analysis and the identification of qualitative 
behaviors. The approach of perturbation methods is 
especially useful because the alternative, running large-scale numerical
simulations from which one may calculate statistics, is particularly
burdensome for slow-fast systems. Such systems generally require the 
use of implicit numerical integrators in order to ensure numerical 
stability, the use of which is extremely computationally expensive.

Separately, stochastic systems are frequently used to model
both microscale and macroscale behaviors that are inherently noisy
or simpler to visualize as driven by randomness. Examples of these
systems range from networks of sensors in noisy environments to the 
control of epidemics. There are many
intricacies under investigation within this field such as finite
noise effects and stochastic resonance \cite{ghjm98} that provide for
much lively research, but will not be our focus in this paper. We 
will in particular study the effect of small noise on the escape 
times for a particle in a multi-scale potential well. To do so, we 
will make use of the variational theory
of large fluctuations as it applies to finding the \emph{most
probable path} along which noise directs a particle to escape
\cite{dykman08}. 

It is well-known that noise has a significant effect on deterministic dynamical
systems.  For example, consider a given initial state in the basin of
attraction for a given attractor, which might be steady, periodic, or
chaotic. Noise can cause the trajectory to cross the deterministic 
basin boundary and move into another, distinct basin of 
attraction~\cite{dyk90,dmsss92,mil96,lumcdy98}. For sufficiently small
noise, basin boundary crossings usually occur near a saddle on the
boundary. However we note that for large noise, such a crossing may 
be determined by the global manifold structure away from the saddle.

This paper will 
consider small noise effects in particular. Specifically, we will 
investigate the effect of arbitrarily small noise on the escape of a 
particle from a potential well. In the small noise limit, one can apply
large fluctuation theory~\cite{feyhib65,dyk90,dmsss92,lumcdy98}; also 
known as large deviation theory used in white noise 
analysis~\cite{FW84,feyhib65,weinanebook}, this approach enables us to 
determine the first passage times in a multi-scale environment. For a
vector field that exhibits dynamics on only one timescale, it is clear how
to use the theory to generate an optimal path of escape. The theory has
been applied to a variety of Hamiltonian and Lagrangian variational 
problems~\cite{wentzell76,Hu1987,Dykman1994d,FW84,GT84,MS93,HTG94} that
do not exhibit singularly perturbed behavior.

For slow-fast systems however, technical issues 
arise while determining the projection of noise restricted to the 
lower dimensional manifold. Several sample based 
approaches have been developed to understand dimension reduction in
systems that have well separated time scales \cite{berglundgentz}. 
The existence of a 
stochastic center
manifold was proven in ~\cite{box89} for systems with certain spectral
requirements. Non-rigorous stochastic normal 
form analyses (which lead to the stochastic center manifold) were 
performed in~\cite{knowie83,nam90,namlin91}. More rigorous 
theoretical treatments of normal form coordinate transformations for 
stochastic center manifold reduction were developed 
in~\cite{arnimk98,arn98,kabanov}. Later, another method of stochastic
normal form reduction was developed in which 
anticipatory convolutions (integrals into the future of the noise 
processes) that appeared in the equations for the slow dynamics were
ignored \cite{rob08}. This latter stochastic normal form technique was possible
because the epidemic model under study permitted certain assumptions 
on the magnitude of the noise projections. The disadvantage of such
assumptions compared to probabilistic methods is that there must be 
guarantees to keep stochastic solutions bounded in the past and future 
\cite{FS2009}, which we may not always have.

We will restrict 
our study to systems with two stable equilibria separated by an 
unstable equilibrium point in phase space; the method of center 
manifold approximations however is not strictly reserved for this 
case. This paper begins by introducing some general theory related to 
slow-fast systems and center manifold reductions. We then review 
large fluctuation theory and how it applies to determining the 
optimal path between invariant manifolds in stochastic systems. Next we
follow many other works and apply the theory to the example of a 
damped Duffing oscillator to compare. 
Finally we compare the switching time estimated via large 
fluctuation theory with numerical results for the example system. We note
that although much work on white noise model reduction is being done using
sample-based methods and asymptotics, our variational approach is more
general in that it may include non-Gaussian noise sources as well.

\section{Theory}
We consider a general $(m+n)$-dimensional dynamical 
system of stochastic differential equations with two 
well-separated timescales and additive noise on the slow 
variables:

\begin{align}
        \label{gen1}
        \bm{\dot{x}} &= \bm{F}(\bm{x},\bm{y}) + \bm{\Phi}(t)\\
        \label{gen2}
        \epsilon \bm{\dot{y}} &= \bm{G}(\bm{x},\bm{y})
\end{align}

\noindent where $\bm{x} \in \mathbb{R}^m$, $\bm{y} \in 
\mathbb{R}^n$; $\bm{\Phi}(t)$ are stochastic terms with 
characteristics depending on the application; 
$\bm{F}:\mathbb{R}^m \times \mathbb{R}^n \rightarrow \mathbb{R}^m$ 
and $\bm{G}:\mathbb{R}^m \times \mathbb{R}^n \rightarrow 
\mathbb{R}^n$ are differentiable functions with equilibrium points at 
the origin, and $\epsilon$ is a small parameter. Such systems are 
known as singularly perturbed or slow-fast systems \cite{guck04} 
with timescales separated by a ratio of $\epsilon$. In this system, 
$\bm{x}$ are slow variables and $\bm{y}$ are fast variables. 
Rescaling $\tau = \epsilon t$ and temporarily removing the stochastic 
terms results in the \emph{layer equations}. Denoting $(\cdot)' = 
\frac{d}{d\tau}$, the deterministic part of Eqs.\ \eqref{gen1}, 
\eqref{gen2} becomes:

\begin{align}
        \label{rescaledF}
        \bm{x}' &= \epsilon \bm{F}(\bm{x},\bm{y})\\
        \label{rescaledG}
        \bm{y}' &= \bm{G}(\bm{x},\bm{y})\\
        \label{rescaledT}
        \epsilon' &= 0.
\end{align}

Note that since $\epsilon$ is treated as a state variable in Eqs.\
\eqref{rescaledF}--\eqref{rescaledT}, then all terms in Eq.\
\eqref{rescaledF} are necessarily nonlinear. If $\bm{G}(\bm{x},\bm{y})$
has a linear part with nonzero determinant, then there exists an 
$m$-dimensional center manifold tangent to the center eigenspace at
the origin. By the implicit function theorem, we may write the 
manifold locally as a function $\bm{h}:\mathbb{R}^m \times \mathbb{R} 
\rightarrow \mathbb{R}^n$: 

\begin{equation}
        \label{cm_h}
        \bm{y} = \bm{h}(\bm{x},\epsilon).
\end{equation}

Following Carr \cite{car81}, the center manifold may
be approximated to arbitrary order by a polynomial series in
$\bm{x}$ and $\epsilon$. All solutions collapse to this manifold at an
exponential rate since it is hyperbolic.

\subsection{Stochastic switching}
\label{sec:stoch_switch}
Stochastic differential equations cannot be described by 
deterministic orbits representing trajectories of a particle through 
phase space. Instead, other approaches are used to qualitiatively 
describe the system. For instance, sample based techniques may 
describe individual realizations in phase space, or families of such 
realizations. Another technique is to find a probability density 
function (pdf) describing the likelihood of finding a particle at a 
given point and time. If the noise is
Gaussian and uncorrelated in time, the dynamics of the pdf 
$\rho(\bm{z},t)$, where $\bm{z} = (\bm{x}; \bm{y})$ is the 
concatenated vector of state variables, are governed by the 
Fokker-Planck equation \cite{gar03}:

\begin{align}
        \nonumber
        \frac{\partial \rho(\bm{z},t)}{\partial t} =& - 
        \sum_{i=1}^{m+n}
        \frac{\partial}{\partial
        z_i} \left(\rho(\bm{z},t) \, \mathcal{F}_i \right)\\
        \label{fp}
        &+ \sum_{i=1}^{m+n} \sum_{j=1}^{m+n} \frac{\partial^2}{\partial
        z_i \, \partial z_j}\left( D_{ij} \, \rho(\bm{z},t)\right)
\end{align}

\noindent where $\bm{\mathcal{F}} = (\bm{F} ; \bm{G})$ is the 
concatenated vector of functions describing the vector field and 
$D_{ij}$ is a diffusion coefficient matrix.

Equation \eqref{fp} relates the time derivative of the probability 
density function $\rho(\bm{z},t)$ with 
expressions involving spatial derivatives of the vector field 
$\bm{\mathcal{F}}$. Presuming
the characterization of the noise and the vector fields are
autonomous, the pdf will asymptotically approach a steady state
distribution that is independent of time. Therefore, we seek steady 
state solutions to Eq.\ \eqref{fp}; that is,

\begin{equation}
        \label{fp_station}
        \sum_{i=1}^{m+n}\frac{\partial}{\partial
        z_i} \left(\rho(\bm{z}) \, \mathcal{F}_i \right) = 
        \sum_{i=1}^{m+n} \sum_{j=1}^{m+n} \frac{\partial^2}{\partial
        z_i \, \partial z_j} \left(D_{ij} \, \rho(\bm{z})\right).
\end{equation}

If the intensity for each noise term is equal and each component is 
uncorrelated, then we may write
$D_{ij} = D \delta_{ij}$. For the system described in Eqs.\ 
\eqref{gen1}, \eqref{gen2}, it is also relevant that $D_{ij} | 
_{i,j>m} = 0$ since additive noise only affects the slow variables; 
this results in

\begin{equation}
        \label{fp_station2}
        \sum_{i=1}^{m+n}\frac{\partial}{\partial
        z_i} \left(\rho(\bm{z}) \, \mathcal{F}_i \right) =
        D \sum_{i=1}^m \frac{\partial^2}{\partial z_i^2} \rho(\bm{z}).
\end{equation}

We will now assume a certain form for the pdf that will allow us to 
solve Eq.\ \eqref{fp_station2} keeping in mind that the goal is to 
analyze stochastically-induced switching. In the small-noise limit, 
transitions between attractors happen only rarely. Therefore, noise 
leading to a transition is considered to be in the tail of the 
probability distribution that governs the amplitude of the noise. A 
stochastically-induced switch is \emph{most likely} to occur in the presence 
of a hypothesized ``optimal noise,'' which has a finite 
likelihood of occurrence. The path that the system travels through
phase space under the influence of the optimal noise is known as the 
``optimal path.'' Such an event follows an exponential 
distribution which we will use as an ansatz to solve Eq.\ 
\eqref{fp_station2}. The WKB ansatz states that $\rho(\bm{z}) \propto \exp 
\left(-\frac{1}{2D} R(\bm{z})\right)$. Applying this to the 
steady-state Fokker-Planck equation \eqref{fp_station2} yields the 
differential equation

\begin{equation*}
        \sum_{i=1}^{m+n} - \frac{\partial R}{\partial z_i} 
        \mathcal{F}_i + 2 D \frac{\partial \mathcal{F}_i}{\partial 
        z_i} = \sum_{i=1}^m -D \frac{\partial R}{\partial z_i^2} + 
        \frac{1}{2}\left( \frac{\partial R}{\partial z_i}\right)^2.
\end{equation*}

Since we are operating in the small-noise limit, any terms multiplied 
by $D$ will be small; for now, we will neglect them leaving the
first order nonlinear equation:

\begin{equation}
        \sum_{i=1}^{m+n} - \frac{\partial R}{\partial z_i} 
        \mathcal{F}_i = \frac{1}{2} \sum_{i=1}^m 
        \left(\frac{\partial R}{\partial z_i}\right)^2.
        \label{action_relation}
\end{equation}

In some cases, solving for $R$ in Eq.\ \eqref{action_relation} is
possible and would result in a stationary pdf for Eqs.\
\eqref{gen1}, \eqref{gen2}. However, combined with the results in
the following section, we will demonstrate that not only is there
a straightforward way to tackle solving for $R$ but also that it
is intimately related with the principle of least action and the
formulation of an optimal path.

\subsection{Formulation of the Optimal Path}

We wish to study the transition rates due to stochastic fluctuation 
between two energy minima. Consider a system with two stable
equilibrium points $\bm{z}_1$ and $\bm{z}_2$ with a saddle point
$\bm{z}_s$ separating them. Since the noise intensity $D$ is small, 
we assume that switching between the two states will be considered a 
``rare event.'' The frequency of such an event is approximately 
determined by the most likely noise to bring the system from 
$\bm{z}_1$ to $\bm{z}_s$, i.e.\ the optimal noise. A realization of 
the optimal noise is calculated to guide the particle to the saddle 
point $\bm{z}_s$, which corresponds to the mean first passage time (MFPT).
The method to calculate this path makes use of Hamiltonian's 
principle. One may predict the switching rate by first 
finding the optimal path between the two states in an expanded phase 
space which accounts for the noise and then calculating the 
dynamical quantity known as the action along that path. For a 
rigorous explanation of this procedure, see \cite{Dykman1994d,weiss}.

The optimal path is the path that is of
minimal action. We write the action of the
noise on the system \eqref{gen1}, \eqref{gen2} as:

\begin{align}
        \nonumber 
        \mathcal{R} 
        [\bm{x},\bm{y},&\bm{\Phi},\bm{\lambda_x},\bm{\lambda_y}] = 
        \frac{1}{2} \int \bm{\Phi}(t) \cdot 
        \bm{\Phi}(t) dt\\
        \nonumber
        &+ \int \bm{\lambda_x} \cdot
        (\dot{\bm{x}} - \bm{F}(\bm{x},\bm{y}) - \bm{\Phi}(t)) dt\\
        \label{action} &+ \int 
        \bm{\lambda_y} \cdot (\epsilon \dot{\bm{y}} - 
        \bm{G}(\bm{x},\bm{y})) dt.
\end{align}

The action integral Eq.\ \eqref{action} represents the total
effect of noise on the system subject to the constraints of the
vector field. The first term involving the action of the noise is
derived by taking a WKB approximation of the Chapman-Kolmogorov 
equation \cite{dykman08} for the infinitesimal noise events
along the path for white noise. The $\bm{\lambda}$ factors are 
Lagrange multipliers, and the terms multiplying them are the 
constraint equations. The integral is calculated along the path for all
time. We note that Eq.\ \eqref{action} is a natural way to describe the
effects of noise from both white and non-Gaussian sources.

To find the functions that minimize the action, we take the first 
variation of the above equation with respect to the 
independent variables and set them equal to zero. This will give five 
sets of equations that when solved will extremize the action 
$\mathcal{R}$. An example of these variational calculations (with 
variation $\xi \in C^1$ bounded) on the action with respect to the 
functions $x_i$ is:

\begin{align}
        \nonumber \frac{\delta \mathcal{R}}{\delta x_i} =& \int 
        \lambda_{x_j}\left(\dot{\xi} -
        \xi \frac{\partial F_j}{\partial x_i}\right) dt
        + \int \lambda_{y_j} \left(-\xi \frac{\partial G_j}{\partial 
        x_i}\right) dt\\
        \label{drdxint} =& \int \xi\left(-\dot{\lambda}_{x_i} -
        \lambda_{x_j} \frac{\partial F_j}{\partial x_i} - \lambda_{y_j} 
        \frac{\partial G_j}{\partial x_i}\right)dt = 0.
\end{align}

Arriving at the second equality involves integrating by parts; since 
the functional derivative restricts the variations $\xi$ to be
bounded, the first term arising from integration by parts vanishes. 
Given Eq.\ \eqref{drdxint}, we have that the function multiplying 
$\xi$ in the integrand must vanish; this yields the differential 
equation:

\begin{equation}
        \label{drdx}
        \dot{\lambda}_{x_i} + \lambda_{x_j} 
        \frac{\partial F_j}{\partial x_i} + \lambda_{y_j} 
        \frac{\partial G_j}{\partial x_i} = 0.
\end{equation}

In the same way, the following equations were derived for the first 
variation with respect to $y_i$, $\lambda_{x_i}$, $\lambda_{y_i}$ and 
$\Phi_i$:

\begin{align}
        \label{drdy} \frac{\delta \mathcal{R}}{\delta y_i} = 0 &\implies
        \epsilon \dot{\lambda}_{y_i} + \lambda_{y_j} \frac{\partial 
        G_j}{\partial y_i} +
        \lambda_{x_j} \frac{\partial F_j}{\partial y_i} = 0\\
        \label{drdl2} \frac{\delta \mathcal{R}}{\delta \lambda_{y_i}} 
        = 0 
        &\implies \epsilon \dot{y}_i = G_i\\
        \label{drdl1} \frac{\delta \mathcal{R}}{\delta \lambda_{x_i}} 
        = 0
        &\implies \dot{x}_i = F_i + \Phi_i\\
        \label{drdeta} \frac{\delta \mathcal{R}}{\delta \Phi_i} = 0
        &\implies \Phi_i = \lambda_{x_i}
\end{align}

\noindent where $i = 1, \dots, m$ and $i = 1, \dots, n$ for the slow 
and fast variables and their conjugate momenta respectively.

To make a connection with Section \ref{sec:stoch_switch}, we will for 
a moment consider the singular limit as $\epsilon \rightarrow 0$ of 
the vector field in Eqs.\ \eqref{drdy}--\eqref{drdeta}. This 
approximation describes the behavior of a particle in the $x_i$ and 
$\lambda_{x_i}$ coordinates after fast transients have died out and 
yields a system known as the ``slow equations.'' The slow equations 
are:

\begin{align}
        \label{slow1}\dot{x}_i &= F_i + \lambda_{x_i}\\
        \label{slow2}\dot{\lambda}_{x_i} &= -\frac{\partial
        F_j}{\partial x_i} \lambda_{x_j} .
\end{align}

The slow equations represent a conservative system. To calculate the 
corresponding Hamiltonian, we note that:

\begin{align*}
        \dot{x}_i &= \frac{\partial \mathcal{H}}{\partial 
        \lambda_{x_i}}\\
        \dot{\lambda}_{x_i} &= -\frac{\partial \mathcal{H}}{\partial 
        x_i}
\end{align*}

where the Hamiltonian is:

\begin{equation}
        \mathcal{H} = F_i \lambda_{x_i} + \frac{1}{2} \lambda_{x_i} 
        \lambda_{x_i}. \label{gen_hamil} 
\end{equation}

Setting $\mathcal{H} = 0$ in Eq.\ \eqref{gen_hamil} verifies an
intriguing relationship: if one identifies $\lambda_{x_i}(\bm{x}) = 
\frac{\partial R(\bm{x})}{\partial x_i}$, Eq.\ \eqref{gen_hamil} and 
Eq.\ \eqref{action_relation} are equivalent for the singular case. 
This confirms our earlier analysis using variational calculus and 
verifies that $R(\bm{z})$ in the Eikonal approximation is indeed the 
action.

The probability of a rare
event occurring described by that approximation is directly 
proportional to the switching rate, or its inverse, the mean first
passage time. This quantity, denoted $T_S$ is inversely
proportional to the switching rate. Since the action will be
calculated along the optimal path, $R = \min \mathcal{R}$ and the
relation to the switching time is

\begin{equation}
        T_S = c \exp(R/2 D).
        \label{switching_time}
\end{equation}

Since the switching rate is proportional to the probability of a 
large fluctuation, there is a proportionality constant $c$ that is 
yet to be determined. The calculation of this prefactor is the 
subject of ongoing research \cite{dyk10}, but is not the
focus of the current work.

\section{Application: the damped Duffing oscillator}
To test the method, we consider a prototypical example for a
double-welled potential---the damped Duffing oscillator. A stochastic 
variant of this oscillator is:

\begin{align}
        \label{xdot} \dot{x} &= y + \eta (t)\\
        \label{ydot} \epsilon \dot{y} &= x - x^3 - y
\end{align}

\noindent where $\epsilon$ is a small parameter and $\eta(t)$ is a noise
source. We will consider the case where $\eta(t)$ represents
uncorrelated Gaussian white noise and is defined by

\begin{equation*}
        \label{whitenoise}
                \langle \eta(t) \eta(t') \rangle = 2 D \delta(t - t').
\end{equation*}

The noise intensity, which controls the width of the distribution of 
noise, is represented by $D = \sigma^2/2$ where $\sigma$ is the standard
deviation of the noise.

Applying Eqs.\ \eqref{drdy}-\eqref{drdeta} to this system, the 
variational equations for the damped Duffing oscillator in Eqs.\ 
\eqref{xdot}-\eqref{ydot} are:

\begin{align}
        \label{duffing_l1} \dot{\lambda_1} &= (3 x^2 - 1) \lambda_2\\
        \label{duffing_l2} \epsilon \dot{\lambda_2} &= \lambda_2 -
        \lambda_1\\
        \label{duffing_y} \epsilon \dot{y} &= x - x^3 - y\\
        \label{duffing_x} \dot{x} &= y + \lambda_1
\end{align}

Following the language of Kaper \cite{kaper}, there are two limits 
over which the system in Eqs.\ \eqref{duffing_l1}--\eqref{duffing_x} 
may be studied. The first involves immediately taking the limit as 
$\epsilon \rightarrow 0$ in the equations, while the latter involves a
rescaling of time and will be considered in the following section. The
first limit yields the slow equations; they are:

\begin{align}
        \label{sing_l1} \dot{\lambda_1} = (3 x^2 - 1) \lambda_1\\
        \label{sing_x} \dot{x} = x - x^3 + \lambda_1.
\end{align}

The critical dynamics in Eqs.\ \eqref{sing_l1}, \eqref{sing_x} have 
the equilibria $(x,\lambda_1) = \left\{(\pm 1, 0), (0,0), \left(\pm 
\frac{1}{\sqrt{3}}, \mp\frac{2}{3 \sqrt{3}}\right)\right\}$.
Note that in the absence of noise, there is a path connecting the 
equilibria along the $x$ axis. For nonzero noise, there is a
heteroclinic connection in the $x, \lambda_1$ plane between the two 
states which represents the 
optimal path---the most likely trajectory for switching between 
the basins at $x = \pm 1$ and $x = 0$. For this
system it is possible to solve for this path explicitly using a 
series of transformations. The optimal path for the $x$ coordinate 
given as a solution to Eqs.\ \eqref{sing_l1}, \eqref{sing_x} is

\begin{equation*}
        x(t) = \pm \frac{1}{\sqrt{1 - A\exp(2t)}},
\end{equation*}

\noindent where $A$ is an arbitrary coefficient to be determined by
the initial condition. Due to symmetry, it suffices to study 
switching between either $x = \pm 1$ and $x = 0$; we choose to 
examine switching from $-1$ to $0$, i.e.\ the negative branch of 
$x(t)$. By inspection it is clear that $A < 0$,
otherwise solutions would cease to exist in finite time. In calculating
the action this coefficient is irrelevant.  
Choosing $A = -1$ (implying $x(0) = \frac{1}{2}$) without loss of 
generality results in the optimal path:

\begin{equation}
        \label{x_path} x(t) = - \left(1 + \exp(2t)\right)^{-1/2}.
\end{equation}

Integrating and solving for the arbitrary unknown
functions, the Hamiltonian for the slow system is:

\begin{equation}
        \label{sing_hamil} \mathcal{H} = (x - x^3) \lambda_1 + \lambda_1^2/2.
\end{equation}

By inspection, we find that $\mathcal{H} = 0$ at both the origin and 
$(x,\lambda_1) = (\pm 1,0)$. The equation for the curve connecting 
the two states is easily obtained from Eq.\ \eqref{sing_hamil}:

\begin{equation*}
        \lambda_1(x(t)) = 2 (x(t)^3 - x(t)).
\end{equation*}

One may calculate the action in the singular case by carrying out the
integral $R(x) = \int_{-1}^0 \lambda_1(x) dx$. However, this would
ignore the dependence of the action on the fast variables; to
approximate this influence, we will resort to center manifold
approximations.

\section{Center manifold reduction}

To analyze Eqs.\ \eqref{duffing_l1}-\eqref{duffing_x}, we will apply
center manifold approximations to reduce the number of dimensions in the
system. The system must first be rescaled to be placed in a form that 
is amenable for this process. To obtain the layer equations, we apply 
the scaling $t = \epsilon \tau$:

\begin{align}
        \label{fast_l1} \lambda_1' &= \epsilon (3 x^2 - 1) \lambda_2\\
        \label{fast_l2} \lambda_2' &= \lambda_2 - \lambda_1\\
        \label{fast_y} y' &= x - x^3 - y\\
        \label{fast_x} x' &= \epsilon (y + \lambda_1)\\
        \label{fast_e} \epsilon' &= 0.
\end{align}

\noindent One benefit of Eqs.\ \eqref{fast_l1}--\eqref{fast_e} is 
that it is no longer singular as
$\epsilon$ vanishes. A second benefit is that Eqs.\
\eqref{fast_l2}--\eqref{fast_y} involve terms that are linear in
the state variables (a space which now includes $\epsilon$) and
that all other equations are purely nonlinear. Therefore,
the hypotheses of the center manifold theorem are satisfied and 
center manifold reductions may be applied to Eqs.\
\eqref{fast_l1}--\eqref{fast_e} to reduce the dimensionality of the
system \cite{car81} \cite{gh}. Since the vector field Eqs.\
\eqref{fast_l1}--\eqref{fast_e} is smooth, we may assume:

\begin{align}
        \label{cm_y} y &= h(x,\lambda_1,\epsilon)\\
        \label{cm_l2} \lambda_2 &= k(x,\lambda_1,\epsilon)
\end{align}

where $h$ and $k$ are differentiable functions of the quantities
specified. Applying the definitions in Eqs.\ \eqref{cm_y}, 
\eqref{cm_l2} to Eqs.\ \eqref{fast_l2}--\eqref{fast_y} and 
substituting the vector fields in Eqs.\ \eqref{fast_l1}, 
\eqref{fast_x} when applying the chain rule, we obtain a system of 
two partial differential equations that may be solved for the 
unknown functions that will define the center manifold. These 
equations are known as the \emph{center manifold conditions}. 
Beginning with the condition resulting from Eq.\ \eqref{fast_y}:

\begin{align}
        \label{k_cond} \left( \frac{\partial k}{\partial x} x'
        + \frac{\partial k}{\partial \lambda_1} \lambda_1' \right) &=
        k(x,\lambda_1,\epsilon) - \lambda_1,
        \intertext{also for Eq.\ \eqref{fast_l2},}
        \label{h_cond} \left( \frac{\partial h}{\partial x} x'
        + \frac{\partial h}{\partial \lambda_1} \lambda_1' \right) &= x
        - x^3 - h(x,\lambda_1,\epsilon).
\end{align}

In general, solving the partial differential Eqs.\ 
\eqref{h_cond}, \eqref{k_cond} will be difficult. However, the center 
manifold reduction method next calls for making approximations
for the functions $h$ and $k$ in terms of polynomials of increasingly 
higher order in their dependent variables. Each variable contributes
to the order of a given term; to represent this, one may consider each
variable scaled by a parameter $\alpha$. The series is truncated at an
arbitrarily specified order in $\alpha$. Explicitly, this means:

\begin{align*}
        h(x,\lambda_1,\epsilon) = c_0 + &
        \alpha\left(c_1 x + c_2 \epsilon + c_3 \lambda_1\right)
        + \alpha^2\left(c_4 x^2 + c_5 x \lambda_1 \right. \\
        & \left. + c_6 x \epsilon + c_7
        \lambda_1^2 + c_8 \lambda_1 \epsilon + c_9 \epsilon^2\right) +
        \ldots \\
        k(x,\lambda_1,\epsilon) = d_0 + &
        \alpha\left(d_1 x + d_2 \epsilon + d_3 \lambda_1\right)
        + \alpha^2\left(d_4 x^2 + d_5 x \lambda_1 \right.\\ 
        &+ \left. d_6 x \epsilon + d_7
        \lambda_1^2 + d_8 \lambda_1 \epsilon + d_9 \epsilon^2\right) +
        \ldots.
\end{align*}

The center manifold expressions to fourth order in $\alpha$ and 
ordered by power in $\epsilon$ are:

\begin{align}
        \nonumber
        h =& x - x^3 - (x + \lambda_1 - 4 x^3 - 3 x^2 \lambda_1) \epsilon\\
        \label{h_comp}&+ (2 x + \lambda_1) \epsilon^2 + \mathcal{O}(|\alpha|^5)\\
        \nonumber
        k =& \lambda_1 + \left( -\lambda_1+3 {x}^{2}\lambda_1 \right) \epsilon\\
        \label{k_comp}&+ \left( 2 \lambda_1+6 x{\lambda_1}^{2} \right) 
        {\epsilon}^{2} - 5 \lambda_1 \epsilon^3 + \mathcal{O}(|\alpha|^5).
\end{align}

Note that setting $\epsilon=0$ in the expressions for $h$
and $k$ in Eqs.\ \eqref{h_comp}, \eqref{k_comp} recovers precisely 
the critical dynamics of Eqs.\ 
\eqref{sing_l1}, \eqref{sing_x}. Since the expressions are given in
powers of $\alpha$, they do not represent an accounting of all terms 
that may be present for a given high order in $\epsilon$. However,
after taking the series in $\alpha$ to high enough order, low-order
terms in $\epsilon$ stop appearing and the resulting series is
treated as one in $\epsilon$.

\begin{figure}[ht!]
        \centering
        \newlength\figureheight 
        \newlength\figurewidth 
        \setlength\figureheight{2in} 
        \setlength\figurewidth{2.5in}
        \input{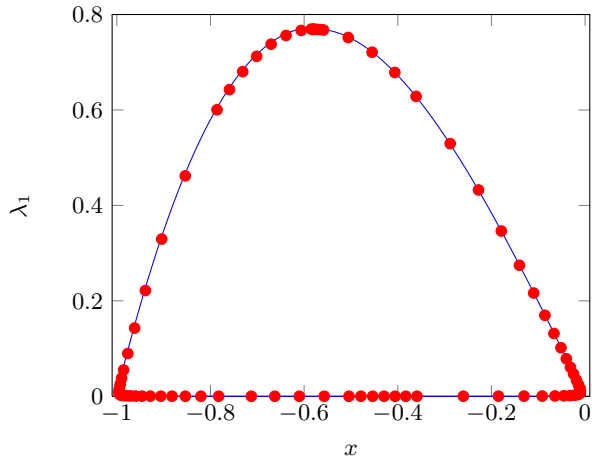}
                \caption{Phase portraits 
                of the full system in Eqs.\ 
                \eqref{duffing_l1}--\eqref{duffing_x} (blue line) 
                compared with its center manifold approximation (red 
                dots) for $\epsilon = 0.001$.}
                \label{fig:duffing_opt_cm_x_px}
\end{figure}

Numerical integration of the original system Eqs.\ 
\eqref{duffing_l1}--\eqref{duffing_x} compared with its center
manifold approximation (with $y$, $\lambda_2$ 
calculated using Eqs.\ \eqref{h_comp}, \eqref{k_comp} respectively) gives 
remarkable agreement even at first order in $\epsilon$. A plot of the 
integration is shown in Figure \ref{fig:duffing_opt_cm_x_px}.

The action of noise along the optimal path in this system is:

\begin{align*}
        \mathcal{R}[x,y,&\eta,\lambda_1,\lambda_2] = 
        \frac{1}{2} \int_{-\infty}^{\infty} \eta^2(t) dt \\
        &+\int_{-\infty}^{\infty} \lambda_1 (\dot{x} - f(x,y) -
        \eta(t)) dt\\
        &+ \int_{-\infty}^{\infty} \lambda_2(\epsilon \dot{y} - 
        g(x,y)) dt
\end{align*}

Having established approximations for these quantities previously
as series expansions in $\epsilon$, it is merely a matter of
careful substitution, differentiation and integration to
obtain an approximation to the integral to arbitrary order in
$\epsilon$. First, substitutions may be made using the center
manifold expressions $y = h(x,\lambda_1,\epsilon)$ and $\lambda_2 = 
k(x,\lambda_1,\epsilon)$ in Eqs.\ \eqref{k_cond},\eqref{h_cond}.
Second, apply the identity along the zero-Hamiltonian curve for 
$\lambda(t)$ as obtained in Eq.\ \eqref{sing_hamil}, along with
$x(t)$ from Eq.\ \eqref{x_path}. Differentiating and integrating as
necessary results in the expression

\begin{equation}
        \label{r0e} \mathcal{R} = \frac{1}{2} - \frac{1}{4}\epsilon^2 
        + \mathcal{O}(\epsilon^3)
\end{equation}

\noindent where the leading order term is the contribution from the 
singular case.

\section{Numerical Results}

To test the predictions resulting from this method, we compared the
scaling predicted from the perturbation method with repeated stochastic
simulation of the damped Duffing oscillator for various values of $D$ 
and $\epsilon$. The stochastic simulations were run using
implicit numerical integration, the details of which are outlined in the
Appendix.

\begin{figure}
        \centering
        \setlength\figureheight{2in} 
        \setlength\figurewidth{2.5in}
%
%
%
%

\definecolor{mycolor1}{rgb}{0.749019622802734,0,0.749019622802734}
\definecolor{mycolor2}{rgb}{0.847058832645416,0.160784319043159,0}

\begin{tikzpicture}

\begin{axis}[%
width=\figurewidth,
height=\figureheight,
scale only axis,
xmin=14,
xmax=29,
xminorticks=true,
xlabel={$1/D$},
ymin=1.5,
ymax=4,
yminorticks=true,
ylabel={$\log_{10} \left(T_S\right)$},
]
\addplot [
color=red,
only marks,
mark=asterisk,
mark options={solid},
forget plot
]
table[row sep=crcr]{
28 3.70100933386968\\
27 3.60946375336655\\
26 3.48993867149049\\
25 3.37132562916638\\
24 3.27952803884116\\
23 3.18226477785982\\
22 3.04718687780978\\
21 2.96202140554888\\
20 2.86388908451343\\
19 2.73223703166328\\
18 2.62940880391491\\
17 2.52728073985499\\
16 2.39606005401962\\
15 2.30781618532135\\
};
\addplot [
color=blue,
only marks,
mark=+,
mark options={solid},
forget plot
]
table[row sep=crcr]{
28 3.6486785150383\\
27 3.55891649741979\\
26 3.44464133941417\\
25 3.35687042314306\\
24 3.24429859743211\\
23 3.12691810884747\\
22 3.04693402477944\\
21 2.92450676059554\\
20 2.79677795986534\\
19 2.67685708386922\\
18 2.58387836074611\\
17 2.47352273401988\\
16 2.36463274383033\\
15 2.27310023360496\\
};
\addplot [
color=mycolor1,
only marks,
mark=x,
mark options={solid},
forget plot
]
table[row sep=crcr]{
28 3.60200804768744\\
27 3.46313156182807\\
26 3.36876393196695\\
25 3.28247113022434\\
24 3.16812779288186\\
23 3.05062339942796\\
22 2.93356867247524\\
21 2.82695563683707\\
20 2.72390046877413\\
19 2.61778607008798\\
18 2.51936811986821\\
17 2.43272680353672\\
16 2.31974014521148\\
15 2.21196627097836\\
};
\addplot [
color=mycolor2,
only marks,
mark=square,
mark options={solid},
forget plot
]
table[row sep=crcr]{
28 3.15304220603509\\
27 3.07865569605835\\
26 2.99054477167532\\
25 2.89258864335514\\
24 2.79020199254775\\
23 2.6780840768839\\
22 2.60879422921215\\
21 2.52781540937489\\
20 2.42231258088697\\
19 2.32043778895142\\
18 2.25290024277462\\
17 2.13942721661401\\
16 2.03387462465458\\
15 1.96357619884927\\
};
\addplot [
color=darkgray,
only marks,
mark=o,
mark options={solid},
forget plot
]
table[row sep=crcr]{
28 2.46874258005738\\
27 2.41135256637961\\
26 2.35581225423176\\
25 2.26651429608434\\
24 2.20517463379458\\
23 2.15570560751883\\
22 2.08098996622498\\
21 2.02465438583125\\
20 1.95661551397988\\
19 1.89560445853698\\
18 1.81985217521188\\
17 1.77982310210292\\
16 1.6778680895528\\
15 1.63400190604065\\
};
\addplot [
color=black,
dashed,
forget plot
]
table[row sep=crcr]{
14 1.64001534333069\\
15 1.6943021535686\\
16 1.7485889638065\\
17 1.80287577404441\\
18 1.85716258428232\\
19 1.91144939452022\\
20 1.96573620475813\\
21 2.02002301499604\\
22 2.07430982523394\\
23 2.12859663547185\\
24 2.18288344570976\\
25 2.23717025594766\\
26 2.29145706618557\\
27 2.34574387642347\\
28 2.40003068666138\\
29 2.45431749689929\\
};
\addplot [
color=mycolor2,
dashed,
forget plot
]
table[row sep=crcr]{
14 1.85002685082871\\
15 1.94502876874504\\
16 2.04003068666138\\
17 2.13503260457772\\
18 2.23003452249405\\
19 2.32503644041039\\
20 2.42003835832673\\
21 2.51504027624306\\
22 2.6100421941594\\
23 2.70504411207574\\
24 2.80004602999207\\
25 2.89504794790841\\
26 2.99004986582474\\
27 3.08505178374108\\
28 3.18005370165742\\
29 3.27505561957375\\
};
\addplot [
color=mycolor1,
dashed,
forget plot
]
table[row sep=crcr]{
14 2.08963007292815\\
15 2.19603222099445\\
16 2.30243436906075\\
17 2.40883651712704\\
18 2.51523866519334\\
19 2.62164081325964\\
20 2.72804296132593\\
21 2.83444510939223\\
22 2.94084725745853\\
23 3.04724940552482\\
24 3.15365155359112\\
25 3.26005370165742\\
26 3.36645584972371\\
27 3.47285799779001\\
28 3.57926014585631\\
29 3.6856622939226\\
};
\addplot [
color=blue,
dashed,
forget plot
]
table[row sep=crcr]{
14 2.15243053322807\\
15 2.26046128560151\\
16 2.36849203797494\\
17 2.47652279034838\\
18 2.58455354272181\\
19 2.69258429509524\\
20 2.80061504746868\\
21 2.90864579984211\\
22 3.01667655221555\\
23 3.12470730458898\\
24 3.23273805696241\\
25 3.34076880933585\\
26 3.44879956170928\\
27 3.55683031408272\\
28 3.66486106645615\\
29 3.77289181882958\\
};
\addplot [
color=red,
dashed,
forget plot
]
table[row sep=crcr]{
14 2.18995468512705\\
15 2.29852287692184\\
16 2.40709106871663\\
17 2.51565926051142\\
18 2.6242274523062\\
19 2.73279564410099\\
20 2.84136383589578\\
21 2.94993202769057\\
22 3.05850021948536\\
23 3.16706841128015\\
24 3.27563660307494\\
25 3.38420479486973\\
26 3.49277298666452\\
27 3.60134117845931\\
28 3.7099093702541\\
29 3.81847756204889\\
};
\end{axis}
\end{tikzpicture}%
        \caption{Mean first passage times from a
        potential well varying with $\epsilon$ and $D$. Data points
        were computed as an ensemble average of 1000 trials. $\circ$ 
        represents $\epsilon = 1.0$, $\color{mycolor2} \Box$ $\epsilon =
        0.5$, $\color{mycolor1} \times$ $\epsilon = 0.2$, $\color{blue} 
        +$ $\epsilon = 0.1$ and $\color{red} \ast$ $\epsilon = 0.01$. 
        Color-corresponding lines show the perturbation-predicted 
        escape times. Lines have been shifted to allow comparison
        with the slopes of simulation data.}
        \label{fig:duffing_escape}
\end{figure}
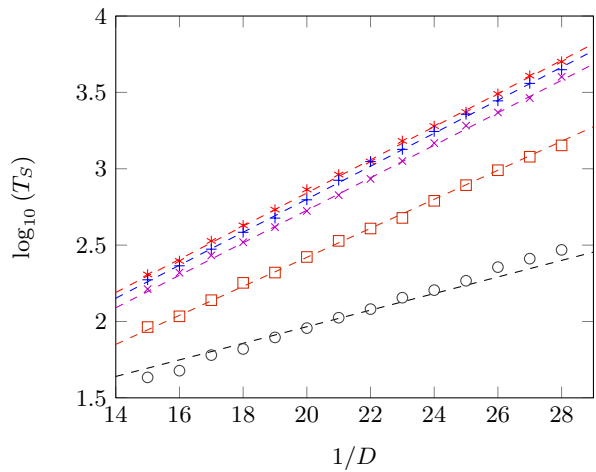

It is convenient to make comparisons between numerics and analytical
approximations by analyzing the logarithm of the escape time across
multiple orders of magnitude. A plot
of the stochastic simulations compared against the escape time as
predicted using the perturbation method is shown in Figure
\ref{fig:duffing_escape}; the two methods agree very
well. Table \ref{table:comparison} provides a side-by-side
comparison of the scaling coefficient between the MFPT and
$\epsilon$ as calculated by the perturbation method and from linear
regression of stochastically simulated switching. The error bounds 
represent the standard deviation on the slope of the regression
line.

\begin{table}[h!]
        \centering
\begin{tabular}{ | c | c | c | }
\hline
$\quad$ & \multicolumn{2}{|c|}{Scaling coefficient $C_S \times
10^{2}$ }\\
\hline
$\epsilon$ & perturbation method & stochastic simulation \\
\hline
\hline
0.001           & 10.86       & 10.91  $\pm$ 1.213   \\
0.003           & 10.86       & 10.84  $\pm$ 1.370   \\
0.01            & 10.86       & 10.79  $\pm$ 1.034   \\
0.1             & 10.80       & 10.80  $\pm$ 1.246   \\
0.2             & 10.64       & 10.60  $\pm$ 1.189   \\
0.5             & 9.500       & 9.295  $\pm$ 1.107   \\
1.0             & 5.428       & 6.469  $\pm$ 0.9437  \\
\hline
\end{tabular}
\caption{Comparison of scaling coefficients for MFPT between the
perturbation method and stochastic simulation. The scaling law is
assumed to be $\log_{10} (T_S) = C_S (1/D) + b$, where
$b$ is a constant determined by simulation. The data above quantify
the predictions and observations in Figure \ref{fig:duffing_escape}.}
\label{table:comparison}
\end{table}

As shown in Table \ref{table:comparison}, the agreement between
stochastic simulation and the perturbation method is quite good and
well within the standard deviation for the slope of the regression
line. However, as the timescales are brought into alignment with one
another, the center manifold approximation applied to the slow
system becomes a poor approximation for the system dynamics. This
can be confirmed visually by observing that the agreement for
$\epsilon = 1.0$ in Figure \ref{fig:duffing_escape} is not
strong.

\section{Discussion}

Figure \ref{fig:duffing_escape} shows that the method fails to predict
the mean first passage time if $\epsilon$ is too large. Both
$\epsilon$ and $D$ are assumed to be small for the perturbation
series and simplifications made. The magnitude of the noise
intensity $D$ may be compared with the height of the barrier
through which the particle must traverse to switch states.
The approximations made do not apply to events where noise is so
significant as to typically cause a transition or where there is little
separation between the time scales. However, the process may be applied to
even higher-dimensional systems where the time scale separation translates
into a spectral gap in relaxation times.

In the regime where $D$ is large compared to the height of the
barrier, the mean first passage time will rapidly decrease. This 
behavior cannot be captured by the WKB approximation ansatz; the 
Eikonal approximation can only capture a linear relationship
between $\log T_S$ and $1/D$. A method could be developed to obtain 
statistics about slow-fast stochastic systems when $D$ is significant 
compared to the effective barrier height, and this will be left to 
future work in which noise is finite and large.

These restrictions aside, the method is resilient to choices of
vector field. Despite the Duffing oscillator's symmetry, the
method has been applied to another double-welled system with broken
symmetry and has resulted in similarly good agreement. Our test system
was an unsymmetric Duffing-like oscillator with differential equations:

\begin{align}
        \label{duffmod1} \dot{x} &= y + \eta(t)\\
        \label{duffmod2} \epsilon \dot{y} &= x (1 + x) (2 - x) - y.
\end{align}

The system in Eqs.\ \eqref{duffmod1}, \eqref{duffmod2} has two stable 
equilibrium points at $x = -1$ and $x = 2$ separated by a saddle at 
$x = 0$. The method outlined in this paper gives the approximate 
expression for the action of:

\begin{equation*}
        \mathcal{R} = \frac{5}{6} - \frac{13}{12} \epsilon^2 +
        \mathcal{O}\left( \epsilon^3 \right).
\end{equation*}

\begin{figure}
        \centering
        \setlength\figureheight{2in} 
        \setlength\figurewidth{2.5in}
        \input{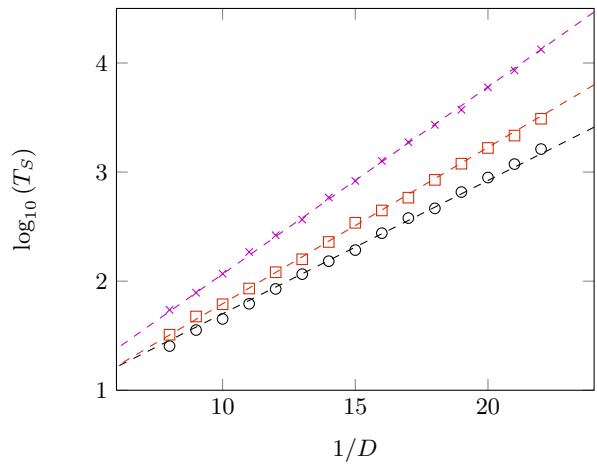}
        \caption{Mean first passage times from a
        potential well varying with $\epsilon$ and $D$. Data points
        were computed as an ensemble average of 1000 trials. $\circ$ 
        represents $\epsilon = 0.5$, $\color{mycolor2} \Box$ 
        $\epsilon = 0.4$, $\color{mycolor1} \times$ $\epsilon = 0.2$.}
        \label{fig:duffmod_escape}
\end{figure}

A comparison of numerically- and formally-generated results 
for the mean first passage time in this system is provided in Figure 
\ref{fig:duffmod_escape}. Both examples we have carried out do not 
have any $\mathcal{O}(\epsilon)$ terms appearing in the approximation to
the action; this may be understood via an analogy with function
optimization. The local behavior of a function at a minimum
with respect to a parameter has no linear dependence on said
parameter by definition.

A contemporary and popular approach to obtain similar results for the
occurrence of rare events uses what are known as sample-based techniques.
Throughout our approach, we have completely avoided the use of such
methods. These approaches generally have required the 
calculation of convolution integrals that depend on the 
realization of noise for all past and future times; while 
analytically tractable, this comes with some assumptions. 
Some of the integrals that result from sample-base approaches must 
remain bounded, putting further restrictions on the noise 
distribution. Such restrictions can be challenging to rigorously 
justify and are at times opaque. Our approach does not require such 
justifications and reaches complementary conclusions while remaining
transparent throughout the process, making it a useful and very
straightforward alternative to sample-based techniques.

Finally, the use of center manifold reductions requires
considerable algebraic manipulation that may not be tenable in all
circumstances, e.g.\ in high dimensional systems or those with
many parameters. Such systems often have lower-dimensional analogs
which may be amenable to this analysis and thus are within reach
of this method. However, there are other approaches. For instance, 
computational methods exist to minimize the action in a variety of
gradient and non-gradient systems \cite{vde04}. These numerical 
algorithms provide a new approach to verify the scaling relationships 
generated by our method in theory and experiment.

\section{Conclusion}

In this work, a method was developed to leverage the disparate
timescales in slow-fast stochastic systems to aid analysis and predict
switching times between attractors. The process avoided the projection
of noise vectors onto the slow manifold in favor of analyzing the 
noisy system via a variational approach to find the optimal path. The 
damped Duffing oscillator was used as an example of a prototypical system
with two potential wells where switching can occur as a result of
large fluctuations. Using this theory, we transformed the original
2-dimensional stochastic system into a 4-dimensional deterministic
system and proceeded to analyze the
optimal path representing the most likely noise to induce a 
transition. The action along this path was crucial to determining
the switching time between the two metastable states present.

For future work, we intend to apply this theory to prescient examples 
of slow-fast stochastic systems, including epidemic models with
non-Gaussian noise. We also will apply this method to systems which 
exhibit delayed feedback.

\section{Acknowledgements}

The authors gratefully thank Luis Mier-y-Teran Romero for helpful
discussions and his prescient insight. This research was performed while
CRH held a National Research Council Research Associateship Award
at the U.S.~Naval Research Laboratory. This research is funded by the 
Office of Naval Research contract F1ATA01098G001 and by Naval Research
Base Program Contract N0001412WX30002.

\appendix
\section{Stochastic simulations}

Ordinary differential equations with multiple timescales present
unique challenges when numerically integrating to obtain a time series,
including the possibility of the system's being ``stiff.'' Stiffness is
a qualitative property of a dynamical system that stymies standard
(i.e.\ explicit) numerical integration methods. This 
effect may be illustrated with a simple example; consider the
system of differential equations:

\begin{align}
        \label{num1} \dot{\bm{x}} &= \bm{F}(\bm{x},\bm{y}) + \alpha 
        \bm{\Phi}\\
        \label{num2} \epsilon \dot{\bm{y}} &= \bm{G}(\bm{x},\bm{y})
\end{align}

\noindent where $\bm{x} \in \mathbb{R}^m$, $\bm{y} \in \mathbb{R}^n$,
$\bm{F}$ and $\bm{G}$ are differentiable functions, $\bm{\Phi}$ is a
white noise term with amplitude controlled by $\alpha$ and $\epsilon$ is a
parameter that tunes the separation of the timescales between the variables
$\bm{x}$ and $\bm{y}$. For the purpose of illustration, we first set $\alpha =
0$. To obtain a time series of Eqs.\ \eqref{num1},
\eqref{num2} there are many numerical recipes that may be applied, 
the simplest of which is Euler's Method. Let $\mathcal{D}$ 
represent taking the Jacobian of the
vector field, and let $\mathcal{D} \bm{F}$, $\mathcal{D} \bm{G}$ be
nonsingular. Euler's Method calls for generating successive 
iterations of the underlying function by discretizing time with a
uniform step size $\nu$ and iterating the resulting map:

\begin{align}
        \label{nummap1} \bm{x}_{k+1} &= \bm{x}_k +
        \nu \bm{F}(\bm{x}_k,\bm{y}_k)\\
        \label{nummap2} \bm{y}_{k+1} &= \bm{y}_k +
        \frac{\nu}{\epsilon} \bm{G}(\bm{x}_k,\bm{y}_k)
\end{align}

In general, the eigenvalues of both $\mathcal{D} \bm{F}$, 
$\mathcal{D} \bm{G}$ are $\mathcal{O}(1)$. Of particular concern is 
the factor of $\frac{\nu}{\epsilon}$, which is generally very large. 
The eigenvalues of Eq.\ \eqref{nummap2} will in general be much 
larger than those of Eq.\ \eqref{nummap1}, which leads to stiffness. 
This inverse relationship between $\nu$ and $\epsilon$ creates a 
numerical quandary since the necessary step size to ensure stability 
is $\mathcal{O}(\epsilon)$, which is arbitrarily small. For accuracy,
step sizes must be chosen much smaller than this necessary step size,
further aggravating the numerical challenges.

To circumvent this complication, implicit methods are often used to
solve for the state of the system after a time step. We now re-introduce
noise by setting $\alpha = 1$ and draw the noise $\bm{W}_k$ at step $k$ 
from a Gaussian distribution with mean 0 and standard deviation 1. Using a 
first-order Milstein method, the implicit recipe used in our 
stochastic simulations is:

\begin{align}
        \label{implicit1} \bm{x}_{k+1} &= (\bm{x}_k,\bm{y}_k) +
        \nu \bm{F}(\bm{x}_{k+1},\bm{y}_{k+1}) + \sqrt{\nu}
        \bm{W}_k\\
        \label{implicit2} \bm{y}_{k+1} &= (\bm{x}_k,\bm{y}_k) +
        \frac{\nu}{\epsilon} \bm{G}(\bm{x}_{k+1},\bm{y}_{k+1})
\end{align}

Solving for $(\bm{x}_{k+1}, \bm{y}_{k+1})$ in Eqs.\ 
\eqref{implicit1}, \eqref{implicit2} is an exercise in nonlinear,
multidimensional root-finding. Since we expect the system's value at two
adjacent timesteps to be close, we may take $\nu$ arbitrarily small such that a
Newton-Raphson iterative scheme will converge to the value of
$\bm{x}_{k+1}$. While the iterative scheme converges quickly, this 
still involves significant computational overhead since it in general 
requires the inversion of a large matrix.

To compute the mean first passage time for the stochastic systems
compared in Fig.\ \ref{fig:duffing_escape}, the implicit recipe provided
in Eqs.\ \eqref{implicit1}, \eqref{implicit2} was used to integrate the
system until noise caused it to escape from the potential well. The
first passage time was recorded, and the system was reset. This ensemble
was run for 1,000 simulations for each value of $\epsilon$ and $D$ and
the mean of all first passage times for the given parameter values was
computed. The total computation time was considerable even on a desktop
PC with 8 processors using \textsc{Matlab}'s parallel computing toolbox.
The processing time to generate Figure \ref{fig:duffing_escape} was over
two weeks using all eight cores clocked at 2.8 GHz.

\bibliography{Refs}

\end{document}